\documentclass[12pt,a4paper]{article}
\usepackage{maa-monthly}

\usepackage{CJKutf8}
\pdfoutput=1

\usepackage{subeqnarray}
\usepackage{multirow}  
\usepackage{tensor}  
\usepackage{float}  

\date{December 31, 2020 \\ revised September 20, 2021}   

\begin{document}

\title{An Elementary Proof of Takagi's Theorem \\
       on the Differential Composition of Polynomials
      }
\markright{Notes}

\author{ Alan D.~Sokal
}

\maketitle

\begin{abstract}
I give a short and completely elementary proof of
Takagi's 1921 theorem on the zeros of a composite polynomial
$f(d/dz) \, g(z)$.
\end{abstract}

%



\newtheorem{theorem}{Theorem}
\newtheorem{proposition}[theorem]{Proposition}
\newtheorem{lemma}[theorem]{Lemma}
\newtheorem{corollary}[theorem]{Corollary}
\newtheorem{definition}[theorem]{Definition}
\newtheorem{conjecture}[theorem]{Conjecture}
\newtheorem{question}[theorem]{Question}
\newtheorem{problem}[theorem]{Problem}
\newtheorem{openproblem}[theorem]{Open Problem}
\newtheorem{example}[theorem]{Example}

\renewcommand{\theenumi}{\alph{enumi}}
\renewcommand{\labelenumi}{(\theenumi)}
\def\eop{\hbox{\kern1pt\vrule height6pt width4pt
depth1pt\kern1pt}\medskip}
\def\prf{\par\noindent{\bf Proof.\enspace}\rm}
\def\rmk{\par\medskip\noindent{\bf Remark\enspace}\rm}

\newcommand{\textbfit}[1]{\textbf{\textit{#1}}}

\newcommand{\bigdash}{%
\smallskip\begin{center} \rule{5cm}{0.1mm} \end{center}\smallskip}

\newcommand{\safepar}{ {\protect\hfill\protect\break\hspace*{5mm}} }

\newcommand{\be}{\begin{equation}}
\newcommand{\ee}{\end{equation}}
\newcommand{\<}{\langle}
\renewcommand{\>}{\rangle}
\newcommand{\widebar}{\overline}
\def\reff#1{(\protect\ref{#1})}
\def\spose#1{\hbox to 0pt{#1\hss}}
\def\ltapprox{\mathrel{\spose{\lower 3pt\hbox{$\mathchar"218$}}
    \raise 2.0pt\hbox{$\mathchar"13C$}}}
\def\gtapprox{\mathrel{\spose{\lower 3pt\hbox{$\mathchar"218$}}
    \raise 2.0pt\hbox{$\mathchar"13E$}}}
\def\textprime{${}^\prime$}
\newcommand{\myendremark}{ $\blacksquare$ \bigskip}
\def\half{ {1 \over 2} }
\def\third{ {1 \over 3} }
\def\twothird{ {2 \over 3} }
\def\smfrac#1#2{{\textstyle{#1\over #2}}}
\def\smhalf{ {\smfrac{1}{2}} }
\newcommand{\real}{\mathop{\rm Re}\nolimits}
\renewcommand{\Re}{\mathop{\rm Re}\nolimits}
\newcommand{\imag}{\mathop{\rm Im}\nolimits}
\renewcommand{\Im}{\mathop{\rm Im}\nolimits}
\newcommand{\sgn}{\mathop{\rm sgn}\nolimits}
\newcommand{\tr}{\mathop{\rm tr}\nolimits}
\newcommand{\supp}{\mathop{\rm supp}\nolimits}
\newcommand{\disc}{\mathop{\rm disc}\nolimits}
\newcommand{\diag}{\mathop{\rm diag}\nolimits}
\newcommand{\tridiag}{\mathop{\rm tridiag}\nolimits}
\newcommand{\AZ}{\mathop{\rm AZ}\nolimits}
\newcommand{\NC}{\mathop{\rm NC}\nolimits}
\newcommand{\PF}{{\rm PF}}
\newcommand{\rk}{\mathop{\rm rk}\nolimits}
\newcommand{\perm}{\mathop{\rm perm}\nolimits}
\def\hboxscript#1{ {\hbox{\scriptsize\em #1}} }
\renewcommand{\emptyset}{\varnothing}
\newcommand{\eqdef}{\stackrel{\rm def}{=}}

\newcommand{\restrict}{\upharpoonright}

\newcommand{\compinv}{{\langle -1 \rangle}}   

\newcommand{\scra}{{\mathcal{A}}}
\newcommand{\scrb}{{\mathcal{B}}}
\newcommand{\scrc}{{\mathcal{C}}}
\newcommand{\scrd}{{\mathcal{D}}}
\newcommand{\scre}{{\mathcal{E}}}
\newcommand{\scrf}{{\mathcal{F}}}
\newcommand{\scrg}{{\mathcal{G}}}
\newcommand{\scrh}{{\mathcal{H}}}
\newcommand{\scri}{{\mathcal{I}}}
\newcommand{\scrj}{{\mathcal{J}}}
\newcommand{\scrk}{{\mathcal{K}}}
\newcommand{\scrl}{{\mathcal{L}}}
\newcommand{\scrm}{{\mathcal{M}}}
\newcommand{\scrn}{{\mathcal{N}}}
\newcommand{\scro}{{\mathcal{O}}}
\newcommand\scroo{
  \mathchoice
    {{\scriptstyle\mathcal{O}}}
    {{\scriptstyle\mathcal{O}}}
    {{\scriptscriptstyle\mathcal{O}}}
    {\scalebox{0.6}{$\scriptscriptstyle\mathcal{O}$}}
  }
\newcommand{\scrp}{{\mathcal{P}}}
\newcommand{\scrq}{{\mathcal{Q}}}
\newcommand{\scrr}{{\mathcal{R}}}
\newcommand{\scrs}{{\mathcal{S}}}
\newcommand{\scrt}{{\mathcal{T}}}
\newcommand{\scrv}{{\mathcal{V}}}
\newcommand{\scrw}{{\mathcal{W}}}
\newcommand{\scrz}{{\mathcal{Z}}}

\newcommand{\bfa}{{\mathbf{a}}}
\newcommand{\bfb}{{\mathbf{b}}}
\newcommand{\bfc}{{\mathbf{c}}}
\newcommand{\bfd}{{\mathbf{d}}}
\newcommand{\bfe}{{\mathbf{e}}}
\newcommand{\bfh}{{\mathbf{h}}}
\newcommand{\bfj}{{\mathbf{j}}}
\newcommand{\bfi}{{\mathbf{i}}}
\newcommand{\bfk}{{\mathbf{k}}}
\newcommand{\bfl}{{\mathbf{l}}}
\newcommand{\bfL}{{\mathbf{L}}}
\newcommand{\bfm}{{\mathbf{m}}}
\newcommand{\bfn}{{\mathbf{n}}}
\newcommand{\bfp}{{\mathbf{p}}}
\newcommand{\bfr}{{\mathbf{r}}}
\newcommand{\bfu}{{\mathbf{u}}}
\newcommand{\bfv}{{\mathbf{v}}}
\newcommand{\bfw}{{\mathbf{w}}}
\newcommand{\bfx}{{\mathbf{x}}}
\newcommand{\bfy}{{\mathbf{y}}}
\newcommand{\bfz}{{\mathbf{z}}}
\renewcommand{\k}{{\mathbf{k}}}
\newcommand{\n}{{\mathbf{n}}}
\newcommand{\vv}{{\mathbf{v}}}
\newcommand{\bv}{{\mathbf{v}}}
\newcommand{\w}{{\mathbf{w}}}
\newcommand{\x}{{\mathbf{x}}}
\newcommand{\y}{{\mathbf{y}}}
\newcommand{\cc}{{\mathbf{c}}}
\newcommand{\zero}{{\mathbf{0}}}
\newcommand{\one}{{\mathbf{1}}}
\newcommand{\bmm}{{\mathbf{m}}}

\newcommand{\ahat}{{\widehat{a}}}
\newcommand{\Zhat}{{\widehat{Z}}}

\newcommand{\C}{{\mathbb C}}
\newcommand{\D}{{\mathbb D}}
\newcommand{\Z}{{\mathbb Z}}
\newcommand{\N}{{\mathbb N}}
\newcommand{\Q}{{\mathbb Q}}
\newcommand{\PP}{{\mathbb P}}
\newcommand{\R}{{\mathbb R}}
\newcommand{\RR}{{\mathbb R}}
\newcommand{\E}{{\mathbb E}}

\newcommand{\ba}{{\bm{a}}}
\newcommand{\bahat}{{\widehat{\bm{a}}}}
\newcommand{\sfa}{{{\sf a}}}
\newcommand{\bb}{{\bm{b}}}
\newcommand{\bc}{{\bm{c}}}
\newcommand{\bchat}{{\widehat{\bm{c}}}}
\newcommand{\bd}{{\bm{d}}}
\newcommand{\bee}{{\bm{e}}}
\newcommand{\beh}{{\bm{eh}}}
\newcommand{\bff}{{\bm{f}}}
\newcommand{\bg}{{\bm{g}}}
\newcommand{\bh}{{\bm{h}}}
\newcommand{\bll}{{\bm{\ell}}}
\newcommand{\bp}{{\bm{p}}}
\newcommand{\br}{{\bm{r}}}
\newcommand{\bs}{{\bm{s}}}
\newcommand{\bu}{{\bm{u}}}
\newcommand{\bw}{{\bm{w}}}
\newcommand{\bx}{{\bm{x}}}
\newcommand{\by}{{\bm{y}}}
\newcommand{\bz}{{\bm{z}}}
\newcommand{\bA}{{\bm{A}}}
\newcommand{\bB}{{\bm{B}}}
\newcommand{\bC}{{\bm{C}}}
\newcommand{\bE}{{\bm{E}}}
\newcommand{\bF}{{\bm{F}}}
\newcommand{\bG}{{\bm{G}}}
\newcommand{\bH}{{\bm{H}}}
\newcommand{\bI}{{\bm{I}}}
\newcommand{\bJ}{{\bm{J}}}
\newcommand{\bM}{{\bm{M}}}
\newcommand{\bN}{{\bm{N}}}
\newcommand{\bP}{{\bm{P}}}
\newcommand{\bQ}{{\bm{Q}}}
\newcommand{\bR}{{\bm{R}}}
\newcommand{\bS}{{\bm{S}}}
\newcommand{\bT}{{\bm{T}}}
\newcommand{\bW}{{\bm{W}}}
\newcommand{\bX}{{\bm{X}}}
\newcommand{\bY}{{\bm{Y}}}
\newcommand{\bIB}{{\bm{B}^{\rm irr}}}
\newcommand{\bOB}{{\bm{B}^{\rm ord}}}
\newcommand{\bOS}{{\bm{OS}}}
\newcommand{\bERR}{{\bm{ERR}}}
\newcommand{\bSP}{{\bm{SP}}}
\newcommand{\bMV}{{\bm{MV}}}
\newcommand{\bBM}{{\bm{BM}}}
\newcommand{\balpha}{{\bm{\alpha}}}
\newcommand{\balphapre}{{\bm{\alpha}^{\rm pre}}}
\newcommand{\bbeta}{{\bm{\beta}}}
\newcommand{\bgamma}{{\bm{\gamma}}}
\newcommand{\bdelta}{{\bm{\delta}}}
\newcommand{\bkappa}{{\bm{\kappa}}}
\newcommand{\bmu}{{\bm{\mu}}}
\newcommand{\bomega}{{\bm{\omega}}}
\newcommand{\bsigma}{{\bm{\sigma}}}
\newcommand{\btau}{{\bm{\tau}}}
\newcommand{\bphi}{{\bm{\phi}}}
\newcommand{\bpsi}{{\bm{\psi}}}
\newcommand{\bzeta}{{\bm{\zeta}}}
\newcommand{\bone}{{\bm{1}}}
\newcommand{\bzero}{{\bm{0}}}

\newcommand{\Cbar}{{\overline{C}}}
\newcommand{\Dbar}{{\overline{D}}}
\newcommand{\dbar}{{\overline{d}}}
\newcommand{\Hbar}{{\overline{H}}}
\newcommand{\Sbar}{{\overline{S}}}


\newenvironment{sarray}{
             \textfont0=\scriptfont0
             \scriptfont0=\scriptscriptfont0
             \textfont1=\scriptfont1
             \scriptfont1=\scriptscriptfont1
             \textfont2=\scriptfont2
             \scriptfont2=\scriptscriptfont2
             \textfont3=\scriptfont3
             \scriptfont3=\scriptscriptfont3
           \renewcommand{\arraystretch}{0.7}
           \begin{array}{l}}{\end{array}}

\newenvironment{scarray}{
             \textfont0=\scriptfont0
             \scriptfont0=\scriptscriptfont0
             \textfont1=\scriptfont1
             \scriptfont1=\scriptscriptfont1
             \textfont2=\scriptfont2
             \scriptfont2=\scriptscriptfont2
             \textfont3=\scriptfont3
             \scriptfont3=\scriptscriptfont3
           \renewcommand{\arraystretch}{0.7}
           \begin{array}{c}}{\end{array}}


\newcommand*\circled[1]{\tikz[baseline=(char.base)]{
  \node[shape=circle,draw,inner sep=1pt] (char) {#1};}}
\newcommand{\ostar}{{\,\textcircled{$\star$}\,}}
\newcommand{\ostarN}{{\,\circledast_{\vphantom{\dot{N}}N}\,}}
\newcommand{\ostarPsi}{{\,\circledast_{\vphantom{\dot{\Psi}}\Psi}\,}}
\newcommand{\starN}{{\,\ast_{\vphantom{\dot{N}}N}\,}}
\newcommand{\starpsi}{{\,\ast_{\vphantom{\dot{\bpsi}}\!\bpsi}\,}}
\newcommand{\starone}{{\,\ast_{\vphantom{\dot{1}}1}\,}}
\newcommand{\startwo}{{\,\ast_{\vphantom{\dot{2}}2}\,}}
\newcommand{\starinfty}{{\,\ast_{\vphantom{\dot{\infty}}\infty}\,}}
\newcommand{\starT}{{\,\ast_{\vphantom{\dot{T}}T}\,}}

\newcommand*{\Scale}[2][4]{\scalebox{#1}{$#2$}}

\newcommand*{\Scaletext}[2][4]{\scalebox{#1}{#2}} 

Many theorems in the analytic theory of polynomials
\cite{Dieudonne_38,Marden_66,Rahman_02,Obrechkoff_03}
are concerned with locating the zeros of composite polynomials.
More specifically, let $f$ and $g$ be polynomials (with complex coefficients)
and let $h$ be a polynomial formed in some way from $f$ and $g$;
under the assumption that the zeros of $f$ (respectively, $g$)
lie in a subset $S$ (respectively, $T$) of the complex plane,
we wish to deduce that the zeros of $h$ lie in some subset $U$.
The theorems are distinguished by the nature of the operation defining $h$,
and the nature of the subsets $S,T,U$ under consideration.

Here we shall be concerned with {\em differential composition}\/:
$h(z) = f(d/dz) \, g(z)$, or $h = f(D) \, g$ for short.
In detail, if $f(z) = \sum\limits_{i=1}^m a_i z^i$
and $g(z) = \sum\limits_{j=1}^n b_j z^j$,
then $h(z) = \sum\limits_{i=1}^m a_i \, g^{(i)}(z)$;
and $D$ denotes the differentiation operator, i.e.\ $Dg = g'$.
The following important result was found by Takagi \cite{Takagi_21}
in 1921, subsuming many earlier results:\footnote{
   See \cite{Honda_75,Iyanaga_90,Iyanaga_01,Kaplan_97,Miyake_07}
   for biographies of Teiji Takagi
\begin{CJK}{UTF8}{ipxm}
   ({高木 貞治}, {\em Takagi Teiji}\/, 1875--1960).
\end{CJK}
%
   Takagi's~papers published in languages other than Japanese
   (namely, English, German, and French)
   have been collected in \cite{Takagi_90}.
}

\begin{theorem}[Takagi]
   \label{thm.differential}
Let $f$ and $g$ be polynomials with complex coefficients,
with $\deg f = m$ and $\deg g = n$.
Let $f$ have an $r$-fold zero at the origin ($0 \le r \le m$),
and let the remaining zeros (with multiplicity) be
$\alpha_1,\ldots,\alpha_{m-r} \neq 0$.
Let $K$ be the convex hull of the zeros of $g$.
Then either $f(D) \, g$ is identically zero,
or its zeros lie in the set
${K + \sum\limits_{i=1}^{m-r} [0,n-r] \alpha_i^{-1}}$.
\end{theorem}

\noindent
Here we have used the notations
$A+B = \{ {a+b} \colon\, {a \in A} \hbox{ and } {b \in B} \}$
and $AB = \{ {ab} \colon\, {a \in A} \hbox{ and } {b \in B} \}$.

Takagi's proof was based on Grace's apolarity theorem \cite{Grace_02},
a fundamental but somewhat enigmatic result in the
analytic theory of polynomials.\footnote{
   For discussion of Grace's apolarity theorem and its equivalents
   --- notably Walsh's coincidence theorem
   and the Schur--Szeg\H{o} composition theorem ---
   see \cite[Chapter~IV]{Marden_66}, \cite[Chapter~VII]{Obrechkoff_03},
   and especially
   \cite[Chapter~3]{Rahman_02}.
}
This proof is also given in the books of
Marden \cite[Section~18]{Marden_66},
Obrechkoff \cite[pp.~135--136]{Obrechkoff_03},
and Rahman and Schmeisser \cite[Sections~5.3 and 5.4]{Rahman_02}.
Here I~give a short and completely elementary proof
of Takagi's theorem.

The key step --- as Takagi \cite{Takagi_21} observed ---
is to understand the case of a degree-1 polynomial $f(z) = z - \alpha$:

\begin{proposition}[Takagi]
   \label{prop.takagi.degree1}
Let $g$ be a polynomial of degree $n$,
and let $K$ be the convex hull of the zeros of $g$.
Let $\alpha \in \C$, and define $h = g' - \alpha g$.
Then either $h$ is identically zero,
or all the zeros of $h$ are contained in $K$ if $\alpha = 0$,
and in $K + [0,n] \alpha^{-1}$ if $\alpha \neq 0$.
\end{proposition}

The case $\alpha = 0$ is the celebrated theorem of Gauss and Lucas
\cite[Section~6]{Marden_66}, \cite[Chapter~V]{Obrechkoff_03},
\cite[Section~2.1]{Rahman_02},
which is the starting point of the modern analytic theory of polynomials.
My proof for general $\alpha$
will be modeled on Ces\`aro's \cite{Cesaro_1885} 1885 proof
of the Gauss--Lucas theorem \cite[pp.~72--73]{Rahman_02},
with a slight twist to handle the case $\alpha \neq 0$.

\begin{proof}[Proof of Proposition~\ref{prop.takagi.degree1}]
Clearly, $h$ is identically zero if and only if either
(a) $g \equiv 0$ or (b) $g$ is a nonzero constant and $\alpha = 0$.
Moreover, if $g$ is a nonzero constant and $\alpha \neq 0$,
then the zero set of $h$ is empty.
So we can assume that $n \ge 1$.

Let $\beta_1,\ldots,\beta_n$ be the zeros of $g$ (with multiplicity),
so that $g(z) = b_n \prod\limits_{i=1}^n {(z - \beta_i)}$ with $b_n \neq 0$.
If $z \notin K$, then $g(z) \neq 0$, and we can consider
\be
   {h(z) \over g(z)}
   \;=\;
   {g'(z) \,-\, \alpha g(z) \over g(z)}
   \;=\;
   \sum_{i=1}^n {1 \over z - \beta_i}
   \:-\: \alpha
   \;.
 \label{eq.proof.differential}
\ee
If this equals zero, then by taking complex conjugates we obtain
\be
   0
   \;=\;
   \sum_{i=1}^n {1 \over \bar{z} - \bar{\beta}_i}
      \:-\: \bar{\alpha}
   \;=\;
   \sum_{i=1}^n {z - \beta_i \over |z - \beta_i|^2}
      \:-\: \bar{\alpha}
   \;,
\ee
which can be rewritten as
$z \,=\,  \sum\limits_{i=1}^n \lambda_i \beta_i \:+\: \kappa \bar{\alpha}$
where
\be
   \lambda_i
   \;=\;
   {|z - \beta_i|^{-2}
    \over
    \sum\limits_{j=1}^n |z - \beta_j|^{-2}
   }
   \;,
   \qquad
   \kappa
   \;=\;
   {1
    \over
    \sum\limits_{j=1}^n |z - \beta_j|^{-2}
   }
   \;.
\ee
Then $\lambda_i > 0$ and $\sum\limits_{i=1}^n \lambda_i = 1$,
so $\sum\limits_{i=1}^n \lambda_i \beta_i \in K$;
and of course $\kappa > 0$.
Moreover, by the Schwarz inequality we have
\be
   |\alpha|^2
   \;=\;
   \left| \sum_{i=1}^n {1 \over z - \beta_i} \right|^2
   \;\le\;
   n \sum\limits_{i=1}^n |z - \beta_i|^{-2}
   \;=\;
   {n \over \kappa}
   \;,
\ee
so $\kappa \le n |\alpha|^{-2}$.
This implies that $\kappa \bar{\alpha} \in [0,n] \alpha^{-1}$
and hence that $z \in K + [0,n] \alpha^{-1}$.
\end{proof}

We can now handle polynomials $f$ of arbitrary degree
by iterating Proposition~\ref{prop.takagi.degree1}:

\begin{proof}[Proof of Theorem~\ref{thm.differential}]
{}From $f(z) = a_m \Bigl(\prod\limits_{i=1}^{m-r} (z - \alpha_i) \! \Bigr) z^r$
it is easy to see that
$f(D) = a_m \Bigl(\prod\limits_{i=1}^{m-r} (D - \alpha_i) \! \Bigr) D^r$.
We first apply $D^r$ to $g$, yielding a polynomial of degree $n-r$
whose zeros also lie in $K$ (by the Gauss--Lucas theorem);
then we repeatedly apply (in any order) the factors $D - \alpha_i$,
using Proposition~\ref{prop.takagi.degree1}.
\end{proof}

\medskip

\noindent
{\bf Remark.}
When $\alpha = 0$, the zeros of $h = g'$ lie in $K$;
so one might expect that when $\alpha$ is small,
the zeros of $h = g' - \alpha g$ should lie near $K$.
But when $\alpha$ is small and nonzero, the set $K + [0,n] \alpha^{-1}$
arising in Proposition~\ref{prop.takagi.degree1} is in fact very {\em large}\/.
What is going on here?

Here is the answer:
Suppose that $\deg g = n$.
When $\alpha = 0$, the polynomial $h = g'$ has degree $n-1$;
but when $\alpha \neq 0$, the polynomial $h = g' - \alpha g$ has degree $n$.
So, in order to make a proper comparison of their zeros,
we should consider the polynomial $g'$ corresponding to the case $\alpha=0$
as also having a zero ``at infinity.''
This zero then moves to a value of order $\alpha^{-1}$
when $\alpha$ is small and nonzero.

This behavior is easily seen by considering the example of
a quadratic polynomial $g(z) = z^2 - \beta^2$.
Then the zeros of $g' - \alpha g$ are
\begin{subeqnarray}
   z
   & = &
   {1 \pm \sqrt{1 + \alpha^2 \beta^2} \over \alpha}
       \\
   & = &
   - {\beta^2 \over 2} \alpha \,+\, O(\alpha^3) \,,\quad
   2 \alpha^{-1} \,+\, O(\alpha)
   \;.
\end{subeqnarray}
So there really is a zero of order $\alpha^{-1}$,
as Takagi's theorem recognizes.

In the context of Proposition~\ref{prop.takagi.degree1},
one expects that $g' - \alpha g$ has one zero of order $\alpha^{-1}$
and $n-1$ zeros near $K$ (within a distance of order $\alpha$).
More generally, in the context of Theorem~\ref{thm.differential},
one would expect that $h$ has $m-r$ zeros of order $\alpha^{-1}$,
with the remaining zeros near $K$.
It is a very interesting problem
--- and one that is open, as far as I know ---
to find strengthenings of Takagi's theorem
that exhibit these properties.
There is an old result that goes in this direction
\cite[Corollary~18.1]{Marden_66}, \cite[Corollary~5.4.1(ii)]{Rahman_02},
but it is based on a disc $D$ containing the zeros of $g$,
which might in general be much larger than the convex hull $K$ of the zeros.

\medskip

\noindent
{\bf Postscript.}
A few days after finding this proof
of Proposition~\ref{prop.takagi.degree1},
I~discovered that an essentially identical argument
is buried in a 1961 paper of
Shisha and Walsh \cite[pp.~127--128 and 147--148]{Shisha_61}
on the zeros of infrapolynomials.
I was led to the Shisha--Walsh paper by a brief citation
in Marden's book \cite[pp.~87--88, Exercise~11]{Marden_66}.
So the proof given here is not new; but it deserves to be better known.


\begin{acknowledgment}{Acknowledgments.}
This research was supported in part by
U.K.~Engineering and Physical Sciences Research Council grant EP/N025636/1.
\end{acknowledgment}

%

\begin{affil}
Department of Mathematics, University College London,
London WC1E 6BT, UNITED KINGDOM \\
and Department of Physics, New York University, New York, NY 10003, USA \\
sokal@nyu.edu
\end{affil}


\addcontentsline{toc}{section}{Bibliography}
\begin{thebibliography}{99}

\bibitem{Cesaro_1885}  Ces\`aro, E. (1885). Solution de la question 1338.
   {\em Nouvelles Annales de Math\'ematiques}\/ ($3^e$ s\'erie).
   4: 328--330.
   \url{www.numdam.org/article/NAM_1885_3_4__328_0.pdf}

\bibitem{Dieudonne_38}  Dieudonn\'e, J.  (1938).
   {\em La Th\'eorie Analytique des Polyn\^omes d'une Variable
     (\`a Coefficients Quelconques)}\/.
   M\'emorial des Sciences Math\'ematiques, fascicule 93.
   Paris: Gauthier-Villars.
   \url{www.numdam.org/issue/MSM_1938__93__1_0.pdf}


\bibitem{Grace_02}  Grace, J.~H. (1902). The zeros of a polynomial.
   {\em Proc. Cambridge Philos. Soc.}\/ 11: 352--357.

\bibitem{Honda_75}  Honda, K. (1975).
   Teiji Takagi: A biography --- on the 100th anniversary of his birth.
   {\em Comment. Math. Univ. St. Paul.}\/ 24(2): 141--167.
   \url{doi.org/10.14992/00010342}

\bibitem{Iyanaga_90}  Iyanaga, S. (1990). On the life and works of Teiji Takagi,
   in \cite[pp.~354--376]{Takagi_90}.

\bibitem{Iyanaga_01}  Iyanaga, S. (2001).
   Memories of Professor Teiji Takagi.
   In: Miyake, K., ed.
   {\em Class Field Theory --- Its Centenary and Prospect}\/.
   Advanced Studies in Pure Mathematics \#30.
   Tokyo: Mathematical Society of Japan, pp.~1--11.

\bibitem{Kaplan_97}  Kaplan, P.  (1997).
   Takagi Teiji et la d\'ecouverte de la th\'eorie du corps de classes.
   {\em Ebisu --- \'Etudes Japonaises}\/. 16: 5--11.
   \url{www.persee.fr/doc/ebisu_1340-3656_1997_num_16_1_973}

\bibitem{Marden_66}  Marden, M.  (1966).
   {\em Geometry of Polynomials\/}, 2nd ed.
   Providence, RI: American Mathematical Society.
   (First edition 1949.)

\bibitem{Miyake_07}  Miyake, K.  (2007).
   Teiji Takagi, founder of the Japanese school of modern mathematics.
   {\em Japanese J. Math.}\/ 2(1): 151--164.

\bibitem{Obrechkoff_03}  Obrechkoff, N.  (2003).
   {\em Zeros of Polynomials}\/.
   Sofia: Marin Drinov Academic Publishing House.
   (Originally published in Bulgarian:
   Obre\v{s}kov, N. (1963). {\em Nuli na Polinomite}\/.
   Sofia: Izdat. B\v{u}lgar. Akad. Nauk.)

\bibitem{Rahman_02}  Rahman, Q.~I., Schmeisser, G.  (2002).
   {\em Analytic Theory of Polynomials}\/.
   Oxford: Clarendon Press.

\bibitem{Shisha_61}  Shisha, O., Walsh, J.~L.  (1961). 
   The zeros of infrapolynomials with some prescribed coefficients.
   {\em J. Analyse Math.}\/ 9: 111--160.

\bibitem{Takagi_21}  Takagi, T.  (1921). Note on the algebraic equations.
   {\em Proc. Phys.-Math. Soc. Japan}\/. 3(11): 175--179.
   \url{doi.org/10.11429/ppmsj1919.3.11_175}
   \ (Reprinted in \cite[pp.~175--178]{Takagi_90}.)

\bibitem{Takagi_90}  Takagi, T.  (1990).
   {\em Collected Papers}\/, 2nd ed.
   Edited and with a preface by S.~Iyanaga, K.~Iwasawa, K.~Kodaira, and
   K.~Yosida.
   Tokyo: Springer-Verlag.
   (Reprinted by Springer-Verlag, Heidelberg, 2014.)

\end{thebibliography}
\end{document}